\documentclass{amsart}

\usepackage[mathscr]{eucal}
\usepackage[all]{xy}
\usepackage{amssymb}
\usepackage{amsthm}
\usepackage{amsmath}

\newtheorem{theorem}{Theorem}
\newtheorem{corollary}{Corollary}
\newtheorem{lemma}{Lemma}
\newtheorem{proposition}{Proposition}

\theoremstyle{remark}

\newcommand{\T}{\mathbb{T}}
\newcommand{\N}{\mathbb{N}}

\newcommand{\R}{\mathbb{R}}
\newcommand{\C}{\mathbb{C}}

\newcommand{\waysubset}{\subseteq\!\subseteq}

\DeclareMathOperator{\im}{im}

\title{Cuntz semigroups of ideals and quotients and a generalized Kasparov Stabilization Theorem}
\author{Alin Ciuperca\and Leonel Robert \and Luis Santiago}

\date{}
\address{Alin Ciuperca: Department of Mathematics, University of Toronto, Toronto, ON, M5S 2E4, Canada.}
\email{ciuperca@math.toronto.edu}

\address{Leonel Robert: Fields Institute, 222 College St., Toronto, ON,  M5T 3J1, Canada.}
\email{lrobert@math.toronto.edu}

\address{Luis Santiago: Department of Mathematics, University of Toronto, Toronto, ON, M5S 2E4, Canada.}
\email{santiago@math.toronto.edu}
\begin{document}

\begin{abstract}
Let $A$ be a C*-algebra and $I$ a closed two-sided ideal of $A$.
We use the Hilbert C*-modules picture of the Cuntz semigroup to investigate the relations between the
Cuntz semigroups of $I$, $A$ and $A/I$. We obtain a relation on two elements of the Cuntz semigroup of $A$ 
that characterizes when they are equal in the Cuntz semigroup of $A/I$.
As a corollary, we show that the Cuntz semigroup functor is exact. 
Replacing the Cuntz equivalence relation of Hilbert modules by their isomorphism, we obtain a generalization
of Kasparov's Stabilization theorem. 
\end{abstract}


\thanks{Leonel Robert was supported by NSERC}

\maketitle

\section{Introduction}
In recent years the Cuntz semigroup has emerged as a powerful invariant
in the classification of C*-algebras, simple and nonsimple (e.g.,  \cite{bpt}, \cite{ciuperca-elliott}, 
\cite{rordam}, \cite{toms}). In \cite{toms}
Andrew Toms provides  examples of simple $AH$ C*-algebras that cannot be distinguished
 by their standard Elliott invariant ($K$-theory and traces) but that have different 
Cuntz semigroups. The first author and G. A. Elliott show in \cite{ciuperca-elliott} that in the nonsimple case, the Cuntz semigroup
is a classifying invariant for all AI C*-algebras (their approach relies on Thomsen's 
classification of AI C*-algebras; see \cite{thomsen}). 

Here we define the Cuntz semigroups, stabilized and unstabilized,
in terms of 
countably generated Hilbert C*-modules
over the algebra, following 
the approach introduced by  K. Coward, G. Elliott and C. Ivanescu in \cite{kgc}. This construction of the Cuntz semigroup is analogous to
the description of $K_0$ in terms of finitely generated projective modules, and is based on an appropriate 
translation of the notion of Cuntz equivalence of positive elements to the context of Hilbert C*-modules. 
Our  investigation is initially motivated by the following question: 
is the Cuntz semigroup of a quotient of a C*-algebra implicitly determined by
the Cuntz semigroup of the algebra? We deduce a satisfactory answer from the inequality in Theorem \ref{formula} below, of interest in its own right. 

Given a  countably generated right Hilbert C*-module $M$ over $A$, let us denote by $[M]$ the element
that it defines in $Cu_s(A)$, the stabilized Cuntz semigroup of $A$. We denote by $Cu(A)$ the subsemigroup
of $Cu_s(A)$ consisting of the  elements $[M]$ that satisfy $M\subseteq A^n$ for some $n$. This last semigroup 
can also be described in terms of positive elements of $A$ (and $M_n(A)$), and is often denoted by $W(A)$.

Let $I$ be a $\sigma$-unital ideal of $A$.
Then $MI$ is a countably generated right Hilbert C*-module over $I$. We will see that $[MI]$ only depends on the equivalence
class of $M$. Therefore we write $[M]I:=[MI]$.

\begin{theorem}\label{formula}
Let $I$ be a $\sigma$-unital, closed, two-sided ideal of the C*-algebra $A$ and
let $\pi\colon A\to A/I$ denote the quotient homomorphism.
Let $M$ and $N$ be countably generated right Hilbert C*-modules over $A$.  
Then $Cu_s(\pi)([M])\leq Cu_s(\pi)([N])$ if and only if
\begin{equation*}
[M]+[N]I\leq [N]+[M]I.
\end{equation*}
\end{theorem}

It follows from Theorem \ref{formula} that $Cu_s(\pi)([M])=Cu_s(\pi)([N])$ if and only if 
$[M]+[N]I=[N]+[M]I$. Adding $[l_2(I)]$ on both sides and using Kasparov's stabilization theorem
we get that 
\begin{equation}\label{equivalence}
Cu_s(\pi)([M])=Cu_s(\pi)([N])\Longleftrightarrow [M]+[l_2(I)]=[N]+[l_2(I)].
\end{equation}
We will show that the map
$Cu_s(\pi)\colon Cu_s(A)\to Cu_s(A/I)$ is surjective. We conclude that the restriction of $Cu_s(\pi)$ to $Cu_s(A)+[l_2(I)]$
is an isomorphism onto $Cu_s(A/I)$.

In the case of the unstabilized Cuntz semigroups, the semigroup $Cu(A/I)$ is obtained  
as the quotient of $Cu(A)$ by the equivalence relation: $[M]\sim_I [N]$ if 
$[M]\leq [N]+[C_1]$ and $[N]\leq [M]+[C_2]$ for some $C_1$ and $C_2$, Hilbert C*-modules
over $I$.  Here the assumption  that the ideal $I$ is $\sigma$-unital
is not needed. This result, which we prove, was first obtained by Francesc Perera in an 
unpublished work. It can
also be deduced from \cite[Lemma 4.12]{kirchberg-rordam}.

A suitable notion of exactness of sequences of ordered semigroups can be  defined such that
the isomorphism $Cu_s(\pi)$ between $Cu_s(A)+[l_2(I)]$ and $Cu_s(A/I)$ implies
the exactness in the middle of the sequence
\[
0 \longrightarrow Cu_s(I) \stackrel{Cu_s(\iota)}{\longrightarrow} Cu_s(A)
\stackrel{Cu_s(\pi)}{\longrightarrow} Cu_s(A/I) \longrightarrow 0.  
\] 
In Theorem \ref{cuntzsexact} we will show that this is a short exact sequence of ordered semigroups, with splittings
of the maps $Cu_s(\iota)$ and $Cu_s(\pi)$.

We can express (\ref{equivalence}) more directly as follows: $M/MI$ and $N/NI$ are Cuntz equivalent as $A/I$-Hilbert C*-modules
if and only if $M\oplus l_2(I)$ and $N\oplus l_2(I)$ are also Cuntz equivalent. 
In Section \ref{proofofkasparov} we obtain an improvement of this result,
with isomorphism of Hilbert C*-modules instead of Cuntz equivalence. We prove the following theorem.

\begin{theorem}\label{kasparov-lift}
Let $A$ be a C*-algebra and $I$ a $\sigma$-unital, closed, two-sided ideal of $A$. Let $M$ and $N$ be 
countably generated right Hilbert C*-modules
over $A$ and suppose that $\phi\colon M/MI\to N/NI$ is an isomorphism of $A/I$-Hilbert C*-modules. Then
there is $\Phi\colon M\oplus l_2(I)\to N\oplus l_2(I)$, isomorphism of Hilbert C*-modules, 
that induces $\phi$ after passing to the quotient.
\end{theorem}    

Taking $I=A$ we get Kasparov's Stabilization Theorem (\cite[Theorem 2.1]{kasparov}). 
The module $M\oplus l_2(I)/MI\oplus l_2(I)$ is canonically isomorphic
to $M/MI$. It is using this identification--applied also to $N$--that 
$\Phi$ induces $\phi$. Theorem \ref{kasparov-lift} is proved by an adaptation
of the proof given by Mingo and Phillips in \cite{mingo-phillips} of Kasparov's Stabilization Theorem. 

In the last two sections we apply Theorem \ref{kasparov-lift}
in the context of multiplier algebras and 
we prove an equivariant version of Theorem \ref{kasparov-lift} assuming that the group is compact.

\section{Preliminaries on Hilbert C*-modules}
Let $M$ and $N$ be right Hilbert C*-modules over a C*-algebra $A$. We shall denote by $K(M,N)$ the norm closure
of the space spanned by the $A$-module maps $\theta_{u,v}\colon M\to N$, $\theta_{u,v}(x) :=v\langle u,x\rangle$, 
$u\in M$, $v\in N$. 
We shall denote by $B(M,N)$ the space of adjointable operators from $M$ to $N$. 
If $T\in B(M,N)$, $\ker T$ and $\im T$ will denote the kernel and the image of $T$ respectively. 
When $M=N$ 
the spaces $K(M,N)$ and $B(M,N)$ are C*-algebras that we shall denote by $K(M)$ and $B(M)$ respectively.
The elements of $B(M,N)$ will often be referred to 
simply as operators, while the elements of $K(M,N)$ will be called compact operators. Sometimes we will drop 
the prefix C* and refer to Hilbert C*-modules as Hilbert modules. 
The C*-algebra will always act on the right of the Hilbert C*-modules.

Given a Hilbert C*-module $M$, the Hilbert C*-module $l_2(M)$ is defined as the sequences 
$(m_i)_{i\in \N}$,  $m_i\in M$, with the property that $\sum_i\langle m_i,m_i\rangle$ converges in norm. This module is endowed with the inner
product $\langle (m_i^1),(m_i^2)\rangle :=\sum_i \langle m_i^1,m_i^2\rangle$. 

Let $I$ be a closed two-sided ideal of $A$. By $MI$ we denote
the span of the elements of the form $m \cdot i$, with  $m \in M$, $i \in I$. This set is a closed submodule of $M$ (by Cohen's Theorem, \cite{lance}) consisting of all vectors $z$ of
$M$ for which $\langle z, z \rangle \in I$. 
The quotient $M/MI$ is a right $A/I$-Hilbert C*-module module with inner product 
$\langle x + MI, y + MI \rangle := \langle x, y \rangle + I$. 

The submodule $MI$ is invariant by any operator $T\in B(M)$. More generally, if $T\in B(M,N)$, then
$T(MI)\subseteq NI$. In this way every operator $T$ induces an operator  
$\tilde\pi(T)\in B(M/MI,N/NI)$. 

We say that a Hilbert C*-module $M$ is countably generated if there is a countable set $\{v_i\}_{i=1}^\infty\subset M$ with 
dense span in $M$. 
We will make use of the following two theorems on countably generated Hilbert modules. 
\begin{theorem}\label{hilbert-tietze} (Noncommutative Tietze extension Theorem for Hilbert C*-modules.)
Let $M$ and $N$ be countably generated Hilbert C*-modules and $\phi\in B(M/MI,N/NI)$. Then there
is  $\Phi\in B(M,N)$ that induces $\phi$ in the quotient.
\end{theorem}
\begin{proof}
Let $H=M\oplus N$. We have $H/HI\simeq M/MI\oplus N/NI$.
Using this isomorphism, we define $\psi\colon H/HI\to H/HI$, adjointable operator,  by the matrix
\[
\psi :=
\begin{pmatrix}
0 & \phi^{*}\\
\phi & 0
\end{pmatrix}.
\] 
The homomorphism $\tilde\pi\colon B(H)\to B(H/HI)$ maps $\theta_{u,v}$ to $\theta_{\pi(u),\pi(v)}$
(here $\pi\colon H\to H/HI$ is the quotient map). Thus, $K(H)$ is mapped surjectively onto $K(H/HI)$ by
$\tilde\pi$. 
Since $H$ is countably generated, $K(H)$ is $\sigma$-unital. Thus, 
by the noncommutative Tiezte extension Theorem (\cite[Theorem 2.3.9]{wegge-olsen}), 
$\tilde\pi$ is also surjective. Let
$\Psi\in B(H)$ be a selfadjoint preimage of $\psi$  given by the matrix
\[
\Psi =
\begin{pmatrix}
A & \Phi^*\\
\Phi & B
\end{pmatrix}.
\]
Then the operator $\Phi$ is a lift of $\phi$.
\end{proof}

The following theorem is due to Michael Frank (\cite[Theorem 4.1]{frank}).
\begin{theorem}\label{frank}
Let $M$ and $N$ be Hilbert C*-modules, $M$ countably generated. Let $T\colon M\to N$ be
a module morphism that is  bounded and bounded from below (but not necessarily adjointable).
Then $M$ is isomorphic to $\im T$ as Hilbert C*-modules. 
\end{theorem}

\section{Cuntz semigroups}\label{cuntzsemigroups}
Let us briefly review  the construction of the Cuntz semigroups, stabilized and unstabilized, of
a C*-algebra $A$, in terms of countably generated Hilbert C*-modules over $A$. We refer to \cite{kgc} 
for further details.

Let $M$ be a Hilbert C*-module over $A$.
A submodule $F$ of $M$ is said to be compactly contained in $M$ if there is $T\in K(M)^+$ such that $T$ restricted to $F$ is the identity of $F$. In this case we write $F\waysubset M$. Given two Hilbert C*-modules $M$
and $N$ we say that $M$ is Cuntz smaller than $N$, denoted by $M\preceq N$, if for all $F$, 
$F\waysubset M$, there is $F'$, $F'\waysubset N$, isomorphic to $F$. This relation defines a preorder relation 
on the isomorphism classes 
of Hilbert modules over $A$. Let us say that $M$ is Cuntz equivalent to $N$ if $M\preceq N$
and $N \preceq M$. Let $[M]$ denote the equivalence class of all the modules Cuntz equivalent  to a given 
module $M$. Then the relation $[M]\leq [N]$ if $M\preceq N$ defines an order on the Cuntz equivalence
classes  of right Hilbert modules over $A$. 

Following \cite{kgc}, the stabilized Cuntz semigroup is defined as the ordered set of Cuntz equivalence classes of 
countably generated Hilbert modules over $A$ endowed with the addition law $[M]+[N]:=[M\oplus N]$.  
We denote this ordered semigroup by $Cu_s(A)$. It is shown in \cite{kgc} that $Cu_s(A)=W(A\otimes K)$, where 
$W(A)$ is the Cuntz semigroup of $A$ defined in terms of positive elements of $\cup_{n=1}^\infty M_n(A)$. 
The unstabilized Cuntz semigroup of $A$, denoted by $Cu(A)$, is defined as the subsemigroup of $Cu_s(A)$
formed by the Cuntz equivalence classes $[M]$ of $A$-Hilbert modules $M$ such that $M\subseteq A^n$ for some $n\geq 1$. It is
shown in \cite{kgc} that this ordered semigroup coincides with $W(A)$. Furthermore, we can define functors
$Cu_s(\cdot)$ and $Cu(\cdot)$ from the category of C*-algebras to the category of ordered semigroups.
By choosing a suitable subcategory of the category of ordered semigroups, Coward, Elliott and Ivanescu
were able to show in \cite{kgc} that the functor $Cu_s(\cdot)$ is continuous with respect to inductive limits.

Let $I$ be a $\sigma$-unital closed two-sided ideal of $A$. If $M$ is a countably generated Hilbert module over $A$
then $MI$ is also countably generated.  
Let us see that if $[M]\leq [N]$ then $[MI]\leq [NI]$. Suppose that $F\waysubset MI$. Then  there is 
$F'\waysubset N$ isomorphic to it. Since $F$ and $F'$ are isomorphic and $FI=F$, we must have $F'I=F'$.
So $F'\subseteq NI$. Hence $[F]=[F']\leq [NI]$. Taking supremum over $F\waysubset MI$ we get that $[MI]\leq [NI]$. In particular,
if $M$ and $M'$ are Cuntz equivalent then $MI$ and $M'I$ are also Cuntz equivalent.
This justifies writing $[MI]:=[M]I$. We have seen already that the map $[M]\mapsto [M]I$ 
is order preserving. Since $(M\oplus N)I=MI\oplus NI$, it is also additive. Notice that $M=MI$ (i.e., $M$ is a
Hilbert $I$-module) if and only if $[M]I=[M]$. If $M\subseteq A^n$ then $MI\subseteq A^n$, so the map
$[M]\mapsto [M]I$ sends elements in $Cu(A)$ to elements in $Cu(A)$.  

Let $\iota\colon I\to A$ and $\pi\colon A\to A/I$ denote the inclusion and quotient homomorphisms. 
The morphisms of ordered semigroups $Cu_s(\iota)$ and 
$Cu_s(\pi)$ are given by 
\begin{align*}
Cu_s(\iota)([H]) &:=[H\otimes_\iota A]=[H],\\ 
Cu_s(\pi)([M]) &:=[M\otimes_\pi A/I]=[M/MI].
\end{align*} 
The restrictions of $Cu_s(\iota)$ and $Cu_s(\pi)$ to $Cu(I)$ and $Cu(A)$ respectively, give 
$Cu(\iota)$ and $Cu(\pi)$.

\emph{Proof of Theorem \ref{formula}.}
The hypothesis $Cu_s(\pi)([M])=Cu_s(\pi)([N])$ says that $M/MI$ and $N/NI$ are Cuntz equivalent
$A/I$-Hilbert C*-modules. We will first show that if $M/MI$ is isomorphic to a submodule of  $N/NI$ then we have 
$[M]+[N]I\leq [N]+[M]I$. 

Let $\phi\colon M/MI\to N'/N'I$ be an isomorphism of $M/MI$ into $N'/N'I$,  a submodule of $N/NI$.
Let $C\colon M/MI\to M/MI$ be an arbitrary compact operator with dense range. This operator
exists because $M$ is countably generated. Then
$\phi'=\phi \circ C$ is compact and satisfies that $\im \phi'^*\phi'$ is dense in $M/MI$.
Since $\phi'$ is compact, it is also a compact operator after composing it with the inclusion
of $N'/N'I$ into $N/NI$. Let us consider $\phi'$ as a compact operator having codomain $N/NI$. 
Let $T\colon M\to N$ be a compact operator that lifts $\phi'$. We have a commutative diagram
\begin{align*}
\xymatrix{
M\ar[r]^T\ar[d] & N\ar[d]\\
M/MI\ar[r]^{\phi'} & N/NI.
}
\end{align*}
Since  $\im \phi'^*\phi'$ is dense in $M/MI$, we have that $\im(T^*T)+MI$ is dense in $M$. Let 
$D_1\colon M\to M$ be positive and with $\im D_1$ dense in $MI$. The operator $D_1$ exists because 
$MI$ is countably generated (here we use that $I$ is $\sigma$-unital). Then $T^*T+D_1$ has dense
range in $M$, that is, it is strictly positive. Let $\{F_n\}_{n=1}^\infty$ be an increasing sequence of submodules
of $M$ such that $T^*T+D_1$ is bounded from below on $F_n$ and $\cup_n F_n$ is dense in $M$ 
(e.g., $F_n=\im \phi_n(T^*T+D_1)$, where $\phi_n\in C_0(\R^+)$ has compact support and $\phi_n(t)\uparrow 1$).
Let $G$ be compactly contained in $NI$. We claim that $F_n\oplus G$ is isomorphic to a submodule on 
$N\oplus MI$. By Theorem \ref{frank}, in order to prove this it is enough to find an operator (not necessarily
adjointable)
$\Phi\colon M\oplus NI\to N\oplus MI$ that is bounded from below when restricted to $F_n\oplus G$.
Let us take 
\[
\Phi :=
\begin{pmatrix}
T & -\iota_{NI,N} \\
D_1 & T^*
\end{pmatrix},
\]
where $\iota_{NI,I}$ is the inclusion map of $NI$ in $N$.
In order to show that $\Phi$ is bounded from below it is enough to show that
$\Phi'\Phi$ is bounded from below, where $\Phi'$ is some bounded--possibly 
nonadjointable--operator. 
Let us choose $\Phi'\colon N\oplus MI\to M\oplus NI$ as follows:
\[
\Phi' :=
\begin{pmatrix}
T^* & \iota_{MI,M} \\
-D_2 & T
\end{pmatrix}.
\]
where $D_2\colon N\to N$ has image in $NI$ and is bounded from below on $G$. Then $\Phi'\Phi$ has the form
\[
\Phi'\Phi=
\begin{pmatrix}
T^*T+D_1 & 0 \\
  * & TT^*+D_2
\end{pmatrix}.
\]
To show that the restriction of $\Phi'\Phi$ to $F_n\oplus G$ is bounded from below 
it is enough to show that the operators on the main
diagonal are bounded from below (because the upper right corner is 0). 
This is true by our choice of $F_n$ and $D_2$. 
So $F_n\oplus G$ is isomorphic to a submodule of $N\oplus MI$. Taking 
supremum over $F_n$ and $G$ we get that $[M]+ [NI]\leq [N]+[MI]$.

Now suppose that $M/MI\preceq N/NI$. Let $F\waysubset M$. Then $F/FI\waysubset M/MI$, so $F/FI$ is isomorphic 
to a submodule of $N/NI$. It follows that $[F]+[N]I\leq [N]+[F]I$. Taking supremum 
over all $F$, $F\waysubset M$, we get $[M]+ [NI]\leq [N]+[MI]$. 
\qed

\begin{corollary}
Let $I$ and $J$ be $\sigma$-unital ideals. Suppose that $[M/M(I\cap J)]=[N/N(I\cap J)]$. Then
\[
[M]I+N[J]=[M](I+J)+[N](I\cap  J)=[M](I\cap J)+[N](I+J).
\]
\end{corollary}
\begin{proof}
We have $[M]I+[N]J=[M(I+J)]I+[NJ]=[M(I+J)]+[NJ]I=[M](I+J)+[N](I\cap J)$.
\end{proof}

\begin{corollary}\label{cuntzquotient}
The map $Cu_s(\pi)$ restricted to $Cu_s(A)+[l_2(I)]$ is an isomorphism onto $Cu_s(A/I)$. 
\end{corollary}

\begin{proof} 
As remarked in the introduction, it follows from Theorem \ref{formula} and Kasparov's Stabilization Theorem
that the map $Cu_s(\pi)$ is injective on $Cu_s(A)+[l_2(I)]$.    
$Cu_s(\pi)$ is surjective, since every $A/I$-Hilbert module can be embedded in $l_2(A/I)$, 
and then have its preimage taken by the quotient map $l_2(A)\to l_2(A/I)$. 
$Cu_s(\pi)$ is also surjective restricted to $Cu_s(A)+[l_2(I)]$, since adding $[l_2(I)]$
does not change the image in $Cu_s(A/I)$. Hence, $Cu_s(\pi)$ sends $Cu_s(A)+[l_2(I)]$ isomorphically onto $Cu_s(A/I)$.
\end{proof}

The description of $Cu_s(A/I)$ obtained in Corollary \ref{cuntzquotient} assumes that the
ideal $I$ is $\sigma$-unital. It is possible to obtain $Cu(A/I)$ as a quotient
of $Cu(A)$ by a suitable equivalence relation without assuming that $I$ is $\sigma$-unital.
Since $Cu_s(A)\simeq Cu(A\otimes K)$, this result can also be applied to the stabilized Cuntz
semigroup.  

Recall that $Cu(A)\simeq W(A)$, the latter semigroup defined as equivalence classes
of positive elements on $\cup_n M_n(A)$ (see \cite{rordam}). Given $[a],[b]\in W(A)$
let us say that $[a]\leq _I [b]$ if there are $c\in M_n(I)^+$ for some $n$ such that
$[a]\leq [b]+[c]$. We say that $[a]\!\!\sim_I [b]$ if  $[a]\leq_I [b]$ and $[b]\leq_I [a]$.

\begin{proposition}
The semigroups  $W(A)/\!\!\sim_I$ and $W(A/I)$ are isomorphic.
\end{proposition}
\begin{proof}
Let $\pi\colon A\to A/I$ be as before, the quotient homomorphism.
Let us show that the map $W(\pi)([a])=\pi([a])$ induces an isomorphism after passing to
the quotient by $\sim_I$. Since $\pi$ is surjective, $W(\pi)$ is also surjective. 
It only rests to show that $W(\pi)([a])\leq W(\pi)([b])$ if and only if $[a]\leq_I [b]$.

Let $a$  and $b$ be positive elements in $M_{n}(A)$, such that $\pi(a) \preceq \pi(b)$. 
For all $k$, there is  $d_k \in M_n(A/I)$ such that $\|\pi(a) - d_k\pi(b)d_k^*\| \leq 1/k$. 
By \cite[Lemma 2.2]{kirchberg-rordam2}, there is $d_k' \in M_n(A/I)$  such that $(\pi(a) - 1/k))_+ =
d_k'\pi(b)d_k'^*$. Let $f_k\in M_n(A)$ be such that $\pi(f_k) = d_k'$. We have 
\[
(a - 1/k)_{+} = f_kbf_k^* + i_k \leq f_kbf_k^* + i_k^+,
\] 
for some $i_k^+\in M_n(I)^+$. We get that $[(a - 1/k)_{+}] \leq
[b] + [i_k^+]$. Let $i \in M_n(I)^+$ be an element such that $[i]$ majorizes the sequence 
$[i_k^+]$ (for example, $i = \sum i_k/(2^k\|i_k^+\|)$. Taking supremum over $k$ in 
$[(a - 1/k)_+] \leq [b] + [i]$ we get $[a]\leq [b] + [i]$. Hence $[a] \leq_I [b]$. 
\end{proof}

\section{Exactness of the Cuntz semigroup functor}
Given $S$ and $T$ ordered, abelian semigroups, and $\phi\colon S \to T$ an order preserving  semigroup map, 
let us define
$\mathrm{Ker}(\phi)$  and $\mathrm{Im}(\phi)$ as follows: 
\begin{align*}
\mathrm{Ker}(\phi) & := \{\, (s_{1},s_{2})\in S \times S \mid \phi(s_{1}) \leq
\phi(s_{2}) \,\},\\ 
\mathrm{Im}(\phi) & := \{\, (t_{1}, t_{2}) \in T \times T \mid \exists\,  s_{1}, s_{2} \in S,\, t_{1}
\leq \phi(s_{2}) + t_{2} \,\}.
\end{align*}
We denote by  $\im \phi$ and $\ker \phi$ the image and the kernel of $\phi$ (i.e., the elements mapped
to 0), in the standard sense. 

By a short exact sequence of ordered semigroups we mean one which is exact with respect to the 
two notions of image and kernel defined above.

\begin{theorem}\label{cuntzsexact}
Let $I$ be a $\sigma$-unital ideal of $A$. The short exact sequence 
\[
0 \longrightarrow I \stackrel{\iota}{\longrightarrow} A 
\stackrel{\pi}{\longrightarrow} A/I \longrightarrow 0,
\]
induces split, short exact sequences of ordered abelian semigroups
\begin{align}
0 \longrightarrow Cu_s(I) \stackrel{r}{\leftrightarrows} Cu_s(A)
\stackrel{q}\leftrightarrows Cu_s(A/I) \longrightarrow 0,\label{shortstable}\\
0 \longrightarrow Cu(I) \stackrel{r}{\leftrightarrows} Cu(A)
\longrightarrow Cu(A/I) \longrightarrow 0.\label{shortunstable}  
\end{align}
These sequences are also exact in the standard sense.

The maps $r$ and $q$ are defined as follows: $r([H]):=[HI]$ and $q([M]):=[M']+[l_2(I)]$, 
where $[M']$ is such that $Cu_s(\pi)([M'])=[M]$.
\end{theorem}

\begin{proof}
We have already shown in Corollary \ref{cuntzquotient} that the maps $Cu_s(\pi)$
and $Cu(\pi)$ are surjective. The maps $Cu_s(\iota)$ and $Cu(\iota)$ are injective, since if 
$M$ is Cuntz smaller than $N$ as $I$-modules, then the same holds when they are regarded as $A$-modules.

Let us prove the exactness of the sequence \eqref{shortstable} and note that the same proof works also for the sequence \eqref{shortunstable}. Exactness at $Cu_s(I)$ and $Cu_s(A/I)$ is easily verified.
To check the exactness in the middle  of the sequence \eqref{shortstable} it suffices to prove that $\mathrm{Ker}(Cu_s(\pi))\subseteq\mathrm{Im}(Cu_s(\iota))$, the other inclusion being obvious. The pair $([M],[N])$ belongs to $\mathrm{Ker}(Cu_s(\pi))$ precisely when $Cu_s(\pi)([M])=Cu_s(\pi)([N])$, and this is equivalent by Theorem \ref{formula} with the fact that $[M]+[N]I=[N]+[M]I$. This shows that $([M],[N])\in\mathrm{Im}(Cu_s(\iota))$, and hence $\mathrm{Ker}(Cu_s(\pi))\subseteq\mathrm{Im}(Cu_s(\iota))$. 

It only remains to show that the maps $q$ and $r$ define splittings
of $Cu_s(\pi)$ and $Cu_s(\iota)$ respectively. We have already observed
that $Cu_s(\pi)$ restricted to $Cu_s(A)+[l_2(I)]$ is an isomorphism of
ordered semigroups. Its inverse is $q$. We have also noted that $M=MI$ (i.e., $M$ is a
Hilbert $I$-module) if and only if $[M]I=[M]$, which shows that $r$ is a splitting of $Cu_s(\iota)$.

The restriction of $r$ to $Cu(A)$ is a splitting of $Cu(\iota)$. 
\end{proof}

\emph{Remarks.} The map $r$ 
does not preserve the way below relation of elements in $Cu_s(A)$ (for the definition of
this relation, see \cite{kgc}). So, it is not a morphism in the
category of ordered semigroups defined by Coward, Elliott and Ivanescu. However, $r$ does preserve directed suprema. 

\begin{proposition}
Let $\{[H_i]\}_{i=1}^\infty$ be an increasing sequence in $Cu_s(A)$ with supremum $[H]$. Then 
$[H]I=\sup_i ([H_i]I)$.
\end{proposition}
\begin{proof}
It will be enough to show that $[H]I\leq \sup_i([H_i]I)$, the other inequality being obvious. Let $F$ be a compactly contained submodule of $HI$. Then $F$ is compactly contained in $H$, hence we conclude that $[F]\leq[H_i]$ for some $i$ (see \cite[Theorem 1]{kgc}). This implies that $[F]\leq[H_i]I$, so $[F]\leq \sup_i ([H_i]I)$. Taking supremum over $F$, we get that
 $[H]I\leq \sup_i([H_i]I)$.
\end{proof}

\section{Proof of Theorem \ref{kasparov-lift}.}\label{proofofkasparov}
\begin{proof}
By Theorem \ref{hilbert-tietze}, there is  an operator 
$T\in B(M,N)$ that lifts $\phi$. The following diagram commutes:
\begin{equation*}
\xymatrix{
M\ar[r]^T\ar[d] & N\ar[d]\\
M/MI\ar[r]^\phi & N/NI.
}
\end{equation*}
The operator $T$ in general will not be an isomorphism. However, by the commutativity of this diagram, and the fact that $\phi^*=\phi^{-1}$, we do have that
\[N=\im T+NI\, \hbox{ and }\,M=\im T^* + MI.\]
We now follow the ideas of Mingo and Phillips's proof of the Stabilization Theorem (\cite[Theorem 1.4]{mingo-phillips}) 
to find $\widetilde T\colon M\oplus l_2(I)\to N\oplus l_2(I)$ such that $\widetilde T$  and $\widetilde T^*$ have dense
range. The desired isomorphism $\Phi$ will be obtained by the polar decomposition of $\widetilde T$.

Since $I$ is $\sigma$-unital, the modules $MI$ and $NI$ are countably generated.
Let $\{\eta_k\}$ and $\{\zeta_k\}$ be infinite sequences of generators of $MI$ and
$NI$ respectively, such that each generator appears infinitely often.
Let us define  operators $\phi_1\colon l_2(I)\to N$, $\phi_2\colon l_2(I)\to l_2(I)$,
and $\phi_3\colon l_2(I)\to M$ by the formulas
\begin{align}
\phi_1((x_k)) := \sum_k \frac{1}{2^k}\eta_k x_k,\quad
\phi_2((x_k)) := (\frac{1}{4^k} x_k ),\quad
\phi_3((x_k)) := \sum_k \frac{1}{2^k} \zeta_k x_k.
\end{align}
Let $\widetilde T\colon M\oplus l_2(I)\to N\oplus l_2(I)$ be defined by the matrix
\[
\widetilde T :=
\begin{pmatrix}
T & \phi_1\\
\phi_3^* & \phi_2
\end{pmatrix}.
\]
Notice that $\widetilde T$ is still a lift of $\phi$.
We have $\widetilde T(0,y)=(\phi_1y,\phi_2 y)$, for $y\in l_2(I)$. It is argued in the proof of 
\cite[Theorem 1.4]{mingo-phillips}, that 
this set is dense in $NI\oplus N$. Thus $NI\oplus N\subseteq \overline{\im \widetilde T}$.
Also, $\widetilde T(x,0)= (Tx,0) +(0,\phi_3^*y)$. So,
$\im T\oplus 0\subseteq \overline{\im\widetilde T}$. We conclude that $\im \widetilde T$ is dense in
$M_2\oplus l_2(I)$. In the same way we show that
$\widetilde T^*$ has dense range. Thus, the operator $\widetilde T$ admits a polar
decomposition of the form $\widetilde T=\Phi(\widetilde T^*\widetilde T)^{1/2}$, with $\Phi$ an isomorphism (see Proposition 15.3.7 \cite{wegge-olsen}).
Passing to the quotients $M/MI$ and $N/NI$, the operator $\widetilde T^*\widetilde T$ induces the identity.
So $\Phi$ lifts $\phi$.  
\end{proof}

\section{Multiplier algebras}
Let $A$ be a $\sigma$-unital algebra and $I$ a $\sigma$-unital closed two-sided ideal of $A$.
In this section we use Theorem \ref{kasparov-lift} to explore the relationship between the multiplier 
algebras $M(A)$ and $M(A/I)$. 

We shall consider $A$ and $I$ as countably generated right Hilbert modules over $A$. We shall identify 
the algebra $K(A)$ with $A$, and the algebra $B(A)$ with $M(A)$.
All throughout this section we make the following two assumptions:

(1) the ideal $I$ is stable,

(2) $A\simeq A\oplus I$ as $A$-Hilbert modules.
 
Let us denote by $s\colon M(A)\to M(I)$ the map given
by restriction of the multipliers of $A$ to the invariant submodule $I$. 
Let $\tilde\pi\colon M(A)\to M(A/I)$ be the extension of the quotient map $\pi\colon A\to A/I$
by strict continuity. Recall that, by the 
noncommutative Tietze extension theorem, 
$\tilde\pi$ is surjective.

Recall the fact that for $p,q\in M(A)\otimes M_n(\C)$ projections, 
the modules $pA^n$ and $qA^n$ are isomorphic if and only if   
$p$ and $q$ are Murray-von Neumann equivalent. 
The following lemma gives an alternative way of expressing conditions (1) and (2) above.

\begin{lemma}\label{piprojections}
The following propositions are equivalent.

(i) The ideal $I$ is a direct summand of $A$ as a right $A$-Hilbert module.

(ii) There is a projection $P_I\in M(A)$ such that $P_IA\subseteq I$ and $s(P_I)$ is Murray-von Neumann 
equivalent to the unit of $M(I)$.

Any two projections of $M(A)$ that satisfy (ii) are Murray-von Neumann equivalent in $M(A)$.
\end{lemma}
\begin{proof}
Suppose we have (i). Let $A=I_1\oplus N_1$, with $I_1\simeq I$ as right Hilbert $A$-modules. 
Let $P_I\in M(A)$ be the projection onto $I_1$. Since $I_1$ is an $I$-module, $I_1I=I_1$, hence 
$P_IA=P_II=I_1\subseteq I$. Since the $I$-module $P_II$ is isomorphic to $I$, it follows
that $P_I$, as an $I$ multiplier, is Murray von Neumann equivalent to the identity of $M(I)$.

Suppose we have (ii). The $I$-modules $I$ and $P_II$ are isomorphic. Hence, they are isomorphic
as $A$-modules. Since $P_IA\subseteq I$, we have $P_II=P_IA$. Hence, $P_IA$ is a direct summand
of $A$ isomorphic to $I$.

If $P_I^{(1)}$ and $P_I^{(2)}$ satisfy (ii) then 
$P_I^{(1)}A=P_I^{(1)}I\simeq I\simeq P_I^{(2)}I=P_I^{(2)}A$. Thus, $P_I^{(1)}$ and $P_I^{(2)}$
are Murray-von Neumann equivalent.
\end{proof}

\begin{proposition}\label{unitary-projections}
Suppose $A$ and $I$ satisfy conditions (1) and (2) above. The following 
propositions are true.

(i) Every unitary of $M(A/I)$ lifts to a unitary of $M(A)$.

(ii) If $p$ and $q$ are projections in $M(A)$ such that 
$\tilde\pi(p)$ and $\tilde\pi(q)$ are Murray-von Neumann equivalent in  $M(A/I)$, then 
$p\oplus P_I$ and $q\oplus P_I$ are Murray-von Neumann equivalent
in $M_2(M(A))$.

(iii) For every projection $p_0\in M(A/I)$ there is $p\in M(A)$ such that $\tilde\pi (p)$
is Murray-von Neumann equivalent to $p_0$. 
\end{proposition}
\begin{proof}
(i) Let $\Phi\colon A\to A\oplus I$ be an $A$-module isomorphism. This map induces
an isomorphism $\phi\colon A/I\to (A\oplus I)/(A\oplus I)I$, and composing with the
canonical identification of $(A\oplus I)/(A\oplus I)I$ and $A/I$, we get a unitary
$\phi'\colon A/I\to A/I$. By Theorem \ref{kasparov-lift}, we can lift this unitary to a unitary 
$\Phi'\colon A\oplus I\to A\oplus I$. Now the map $\Phi_0=(\Phi')^{-1}\Phi$ is an isomorphism
of the Hilbert modules $A$ and $A\oplus I$ that induces the identity in the quotient.
    
Let $u\in M(A/I)$ be unitary. By Theorem \ref{kasparov-lift}, there is a unitary $U\colon A\oplus I\to A\oplus I$
that lifts $u$. Then $\Phi_0^*U\Phi_0\in M(A)$ is a unitary that lifts $u$. 

(ii) Since the $A/I$-modules $\pi(p)A/I$ and $\pi(q)A/I$ are
isomorphic, we get $pA\oplus I\simeq qA\oplus I$.
We have $P_IA\simeq I$. Hence,
 $pA\oplus P_IA\simeq qA\oplus P_IA$. So $p\oplus P_I$ is 
Murray-von Neumann equivalent to $q\oplus P_I$. 

(iii) The  $A$-modules $\pi^{-1}(p_0A/I)\oplus \pi^{-1}((1-p_0)A/I)$
and $A$ are isomorphic in the quotient (to $A/I$). Thus 
$\pi^{-1}(p_0A/I)\oplus \pi^{-1}((1-p_0)A/I)\oplus I\simeq A\oplus I\simeq A$.
So $\pi^{-1}(p_0A/I)$ is a direct summand of $A$. Let $p\in M(A)$ be such
that $pA\simeq \pi^{-1}(p_0A/I)$. Then $\pi(p)A/I\simeq p_0A/I$, so 
$\pi(p)$ is Murray-von Neumann equivalent to $p_0$. 
\end{proof}

\emph{Remarks.} If $A$ and $I$ satisfy conditions (1) and (2), then 
$M_n(A)$ and $M_n(I)$ satisfy them as well. So Proposition \ref{unitary-projections}  applies to the pair 
$(M_n(A)$, $M_n(I))$.
If $A$ is stable then $A\oplus I\simeq A$ (by the
Stabilization Theorem), and $I$ is stable. So (1) and (2) are verified in this case too.
More generally, suppose there is $B$ stable such that $I\subseteq A\subseteq B$, and $I$ is an ideal
of $B$. Then there is $P_I\in M(B)$ that satisfies (ii) of Lemma \ref{piprojections}. The restriction
of $P_I$ to $A$ is in $M(A)$ and satisfies (ii) of Lemma \ref{piprojections}. Hence, in this case
the pair $A$, $I$ satisfies (1) and (2).

Proposition \ref{unitary-projections} has K-theoretical implications. 
Part (1),  applied to the algebras $M_n(A)$, implies that the map $K_1(M(A))\to K_1(M(A/I))$ is surjective. 
Parts (ii) and (iii) imply that the map $K_0(M(A))\to K_0(M(A/I))$ is an isomorphism. 
We can improve these results as follows.

Let $B$ be a unital C*-algebra. Let $A\otimes B$ be the minimal tensor product of $A$ and $B$.
Given $H$ and $E$, Hilbert modules over $A$ and $B$ respectively, let us denote 
by $H\otimes E$ the external tensor product of $H$ and $E$ (see \cite{lance}). 
This is an $A\otimes B$ Hilbert module. Given $A$-Hilbert modules $H_1$ and $H_2$, $B(H_1,H_2)\otimes B$ denotes the norm closed subspace of 
$B(H_1\otimes B,H_2\otimes B)$ generated by operators of the form $T\otimes b$, with $T\in B(H_1,H_2)$
and $b\in B$. Note that the composition of operators in $B(H_1\otimes B,H_2\otimes B)$ with operators in 
$B(H_2\otimes B,H_3\otimes B)$ results in operators in $B(H_1\otimes B,H_3\otimes B)$.

Let $M(A,I)$ be the kernel of $\tilde\pi\colon M(A)\to M(A/I)$. We have 
$M(A,I)=\{\, x\in M(A)\mid xa,ax\in I \hbox { for all } a\in A \,\}$. 

\begin{proposition}\label{tensorB}
Let $B$ be a unital C*-algebra and $A$ and $I$ as before. 
Let $p\in M(A,I)\otimes B$ be a projection and $P_I'=P_I\otimes 1$, with 
$P_I$ as in Lemma \ref{piprojections} (ii). Then
$p\oplus P_I'$ is Murray-von Neumann equivalent to $0\oplus P_I'$. 
\end{proposition}
\begin{proof}
The multiplier projection $p$ is an operator from $A\otimes B$ to $A\otimes B$
with range contained in $I\otimes B$. Let $\widetilde p\in B(A,I)\otimes B$ denote
the adjointable operator obtained by simply restricting the codomain of $p$ to $I\otimes B$.
Let $\widetilde P_I'\in B(A,I)\otimes B$ be the corresponding operator for $P_I'$. Notice that
$\widetilde p^*\widetilde p=p$, $\widetilde p\widetilde p^*=s(p)\in M(I)\otimes B$, and
similarly for $\widetilde P_I'$. By \cite[Lemma 16.2]{wegge-olsen}, there is  $V\in M_2(M(I)\otimes B)$, partial isometry, 
such that $V^*V=s(p)\oplus s(P_I')$ and $VV^*=0\oplus s(P_I')$. Let $W$ be defined as
\[
W :=
\begin{pmatrix}
0 & 0\\
0 & (\widetilde P_I')^*
\end{pmatrix}
V
\begin{pmatrix}
\widetilde p & 0\\
0 & \widetilde P_I'
\end{pmatrix}.
\] 
Then $W^*W=p\oplus P_I'$, $WW^*=0\oplus P_I'$, and $W\in M(A,I)\otimes B$. 
\end{proof}

\begin{corollary}
We have 
\begin{align*}
K_i(M(A,I)) &=0, \\
K_i(M(A)) &\simeq K_i(M(A/I)), \\
K_{i}(M(A,I)/I) &\simeq K_{1-i}(I), 
\end{align*}
for $i=0,1$.
\end{corollary}
\begin{proof}
From Proposition \ref{tensorB} we deduce that $K_0(M(A,I)$. Taking $B=C(\T)$,
we get $K_1(M(A,I))=0$. Now by the six term exact sequence associated to the extension
$M(A,I)\to M(A)\to M(A/I)$, we have
 $K_i(M(A))\simeq K_i(M(A/I))$, $i=0,1$. Looking at the extension
$I\to M(A,I)\to M(A,I)/I$, we get that $K_i(I)=K_{1-i}(M(A,I)/I)$, $i=0,1$.
\end{proof}

\emph{Question.} If $I=A$ then $A$ is stable, so the unitary group of $M(A)$ is contractible 
by the Kuiper-Mingo Theorem (see \cite[Theorem 16.8]{wegge-olsen}). Is the unitary group of 
$M(A,I)\,\widetilde{\,}$ contractible in the norm or strict topologies?

\section{Equivariant version of Theorem \ref{kasparov-lift}}
Let $G$ be a locally compact (Hausdorff) group acting on the C*-algebra $A$.
A $G-A$ Hilbert C*-module, or simply a $G-A$-module, is a right Hilbert C*-module endowed with a
continuous action of $G$ such that 
\begin{align*}
\langle g\cdot v,g\cdot w\rangle &=g(\langle v,w\rangle),\\
g\cdot (va) &=(g\cdot v)(g(a)), \hbox{ for all }g\in G,\, v\in M, \hbox{ and }a\in A.
\end{align*}
An operator between $G-A$-modules is equivariant if $T(g\cdot v)=gTv$. The action of $g$ on $T$
is defined as $(g\cdot T)(v)=gT(g^{-1}v)$. $T$ is $G$-continuous if 
the map $g\mapsto g\cdot T$ is continuous in the norm of operators.

Given a $G-A$ module $M$ we denote by  $L_2(G,M)$ the Hilbert C*-module $L_2(G)\otimes M$, where
$L_2(G)$ is the left regular representation of $G$. The action of $G$ on $L_2(G,M)$ is defined as 
$g\cdot (\lambda\otimes m)=(g\cdot \lambda \otimes g\cdot m)$. The $G-A$-module
$L_2(G,M)$ can also be viewed as the completion of $C_c(G,M)$--the $M$-valued continuous 
functions on $G$ with compact support--with respect to the $A$-valued inner product $\langle h_1,h_2\rangle=\int \langle h_1(g),h_2(g)\rangle \,dg$.

Let $I$ be a $\sigma$-unital, closed, two-sided ideal of $A$ that is invariant by the action of $G$. Then
we can define a quotient action of $G$ on $A/I$. More generally, given a $G-A$-module $M$, we can define
a natural (quotient) structure of $G-A/I$-module on $M/MI$.

We now state an equivariant version of Theorem \ref{kasparov-lift} for compact groups
(\cite[Theorem 2.1]{kasparov} and \cite[Theorem 2.5]{mingo-phillips} in the case $I=A$). 
\begin{theorem}
Suppose that the group $G$ is compact. Let $I$ be a $\sigma$-unital, invariant, closed, two-sided ideal of $A$. 
Let $M$ and $N$ be countably generated $G-A$ modules. 
Let $\phi\colon M/MI\to N/NI$ be an equivariant isomorphism.
Then there is $\Phi\colon M\oplus L_2(G,l_2(I))\to N\oplus L_2(G,l_2(I))$,
equivariant isomorphism, that induces $\phi$ in the quotient.
\end{theorem}

\begin{proof}
The proof is an adaptation of the proof of Theorem \ref{kasparov-lift}.  
The equivariant isomorphism $\phi\colon M/MI\to N/NI$ can be lifted to an equivariant
operator $T\colon M\to N$ by first lifting it to an arbitrary operator $T'$, and then averaging
over the group: $Tx=\int (g\cdot T')x \,dg$. (This integration is possible because
for all $x\in M$, the function $(g\cdot T')x$ is continuous in $G$.) 

Next we construct the operator $\widetilde T$, this time making sure it is equivariant. 
For this we need to replace the sequences of vectors $\{\eta_k\}$, $\{\zeta_k\}$,
generators of $MI$ and $NI$, by equivariant operators $\eta_k\colon L_2(G,I)\to M$,
$\zeta_k\colon L_2(G,I)\to N$, such that $\sum \im \eta_k$ is dense in $MI$ and 
$\sum \im \zeta_k$ is dense in $NI$. This is guaranteed by the following lemma.

\begin{lemma}
Suppose that $H$ is a countably generated $G-A$ module. Let $I$ be as before. Then there is a sequence
$\eta_k\colon L_2(G,I)\to H$ of $G$-continuous maps such that 
$\sum \im(\eta_k)$ is dense in $HI$. If $G$ is compact these maps can be chosen equivariant.
\end{lemma}

Before proving the lemma, let us proceed with the proof of the theorem. We define the maps $\phi_1$ and $\phi_3$
replacing the vectors $\eta_k$ and $\zeta_k$ for the operators obtained using the lemma. The definition of the 
map $\phi_2$ is unchanged. The resulting operator $\widetilde T$ is equivariant. Following the same argument
of Mingo and Phillips, $\widetilde T$ and $\widetilde T^*$ have dense range. 
Since the unitary part of an equivariant operator is also equivariant, we get the equivariant
isomorphism $\Phi$ by polar decomposition of $\widetilde T$.

Let us prove the lemma. First suppose that $G$ is only locally compact. It is enough to find a 
$G$-continuous operator from $l_2(L_2(G,I))$ to $H$ with range
dense in $HI$. Let $C_1\colon  l_2(L_2(G,I))\to HI$ be a 
$G$-continuous, surjective operator. Its existence is guaranteed by
the Stabilization Theorem. Let $C_2\in K(HI)$ with dense range. Then $C_2C_1\colon l_2(L_2(G,I))\to HI$ has dense 
range, and since it is compact, it is still an adjointable operator after composing it with the inclusion of $HI$ in $H$.
If $G$ is compact we need to choose $C_1$ and $C_2$ equivariant. $C_1$ exists by the Stabilization Theorem. We take 
$C_2=\int (g\cdot C_2')dg$, with $C_2'\in K(HI)^+$ strictly positive. Then $C_2$ is also strictly positive, thus of dense
range.   
\end{proof}

\emph{Remark.} In the case that $G$ is locally compact, Kasparov (\cite{kasparov}), and Mingo and Phillips
(\cite{mingo-phillips}), obtain a $G$-continuous isomorphism of $M\oplus L_2(G,l_2(A))$ and $L_2(G,l_2(A))$.
Thus, it would be desirable to have a $G$-continuous version of Theorem \ref{kasparov-lift}. 
It is possible to obtain a $G$-continuous lift $T$ of $\phi$. Furthermore, the construction of the
operator $\widetilde T$ can be carried through. However, the proof breaks down at the last step, 
since the unitary part of a $G$-continuous operator need not be $G$-continuous.

\end{document}